\documentclass[11pt]{article}
\usepackage[a4paper]{anysize}\marginsize{3.5cm}{3.5cm}{1.3cm}{2cm}
\pdfpagewidth=\paperwidth \pdfpageheight=\paperheight
\usepackage{amsfonts,amssymb,amsthm,amsmath,eucal}
\usepackage{pgf}
\usepackage{bbm,array}
\usepackage{tikz} 
\usepackage{float}
\restylefloat{table}
\usepackage{subfigure}
\usepackage{caption}

\usetikzlibrary{arrows}
\theoremstyle{plain}
\newtheorem{thm}{Theorem}[section]
\newtheorem{theorem}[thm]{Theorem}

\newtheorem{lemma}[thm]{Lemma}

\theoremstyle{definition}

\newtheorem{thevarthm}[thm]{\varthmname}

\newenvironment{varthm*}[1]{\trivlist\item[]{\bf #1.}\it}{\endtrivlist}

\renewcommand\geq{\geqslant}

\renewcommand\leq{\leqslant}

\newcommand\be{\begin{eqnarray*}}
\newcommand\ee{\end{eqnarray*}}

\newcommand\C{\mathbb C}

\newcommand\F{\mathbb F}

\renewcommand\P{\mathbb P}

\newcommand\call{{\mathcal L}}

\newcommand\newop[2]{\def#1{\mathop{\rm #2}\nolimits}}
\newop\log{log}
\newop\ord{ord}
\newop\Gal{Gal}
\newop\SL{SL}
\newop\GL{GL}
\newop\Bl{Bl}
\newop\mult{mult}
\newop\mass{mass}
\newop\div{div}
\newop\codim{codim}
\newop\sing{sing}
\newop\vdim{vdim}
\newop\edim{edim}
\newop\Ass{Ass}
\newop\size{size}
\newop\reg{reg}
\newop\areg{areg}
\newop\asreg{asreg}
\newop\satdeg{satdeg}
\newop\supp{supp}
\newop\gin{gin}
\newop\ini{in}
\newop\vol{vol}
\newop\sat{sat}
\newop\length{length}
\newop\depth{depth}
\newop\characteristic{char}
\newcommand\eqnref[1]{(\ref{#1})}

\newcommand\eps{\varepsilon}

\def\keywordname{{\bfseries Keywords}}%
\def\keywords#1{\par\addvspace\medskipamount{\rightskip=0pt plus1cm
\def\and{\ifhmode\unskip\nobreak\fi\ $\cdot$
}\noindent\keywordname\enspace\ignorespaces#1\par}}
\def\subclassname{{\bfseries Mathematics Subject Classification
(2000)}\enspace}
\def\subclass#1{\par\addvspace\medskipamount{\rightskip=0pt plus1cm
\def\and{\ifhmode\unskip\nobreak\fi\ $\cdot$
}\noindent\subclassname\ignorespaces#1\par}}

\newcommand\rounddown[1]{\left\lfloor#1\right\rfloor}

\definecolor{qqqqff}{rgb}{0,0,0}
\definecolor{uuuuuu}{rgb}{0,0,0}
\definecolor{zzttqq}{rgb}{0,0,0}
\definecolor{xdxdff}{rgb}{0,0,0}

\begin{document}

\author{M.~Dumnicki, \L .~Farnik, A.~G\l \'owka, M.~Lampa-Baczy\'nska, G.~Malara,\\
T.~Szemberg\footnote{TS was partially supported by NCN grant UMO-2011/01/B/ST1/04875}, J.~Szpond, H.~Tutaj-Gasi\'nska}
\title{Line arrangements with the maximal number of triple points}
\date{\today}
\maketitle
\thispagestyle{empty}

\begin{abstract}
   The purpose of this note is to study configurations of lines
   in projective planes over arbitrary fields having the maximal number
   of intersection points where three lines meet. We give precise
   conditions on ground fields $\F$ over which such extremal configurations exist.
   We show that there does not exist a field admitting a configuration of $11$ lines with $17$ triple points,
   even though such a configuration is allowed combinatorially. Finally, we
   present an infinite series of configurations which have a high number of triple
   intersection points.
\subclass{52C30, 05B30, 14Q10}
\end{abstract}


\section{Introduction}
\label{intro}

   Configurations of points and lines have been the classical object of study
   in geometry. They come up constantly in various branches of contemporary
   mathematics, among others serving as a rich source of interesting
   examples and counter-examples. By way of example, in algebraic geometry,
   arrangements of lines have been studied recently by Teitler in \cite{Tei07}
   in the context of multiplier ideals, in \cite{DST13} as counter-examples
   to the containment problem for symbolic powers of ideals of points
   in the complex projective plane and in \cite{BNAL} in the setup
   of the Bounded Negativity Conjecture and Harbourne constants.

   In combinatorics point line arrangements are subject of classical
   interest and current research. Notably, in the last year we have witnessed
   a spectacular proof of a long standing conjecture motivated by
   the Sylvester-Gallai theorem on the number of ordinary lines. It has
   been proved by Green and Tao \cite{GreTao13} with methods closely
   related to the real algebraic geometry. Symmetric $(n_k)$ configurations
   (mostly in the real Euclidean plane)
   are another classical topic of study in combinatorics. Such configurations
   with triple points (i.e. $k=3$) are well understood, see the beautiful
   monograph by Gr\"unbaum \cite{Gru09}. The classification of $(n_4)$
   configurations has witnessed much progress in recent years and is
   almost completed, mainly due to works of Bokowski and his coauthors,
   see e.g. \cite{BokScg13}, \cite{BokPil14} for
   an up to date account on that path of research. In the present note
   we take up a slightly different point of view and investigate what the
   maximal number of triple points in a configuration of $s$ lines
   over an \emph{arbitrary} ground field is. In Theorem \ref{thm: main},
   which is our main result, we give a full classification for $s\leq 11$.
   The interest in this particular bound is explained by the observation
   that it is the first value of $s$ for which a configuration with
   maximal combinatorially possible number of triple points cannot be realized over any field.
   To the best of our knowledge this is the first
   example of this kind. More precisely, non-realizable configurations
   were known previously, probably the first example was found by
   Lauffer \cite{Lau54}, see also \cite{Stu91} and \cite{Gly99}. However, the known
   examples concern configurations whose numerical invariants (for example
   Lauffer's example is numerically the Desargues $(10_3)$ configuration)
   allow several combinatorial realizations and \emph{some} of them are not
   realizable over any field. In the case studied here, there are two
   combinatorial realizations possible but \emph{none} of them has
   a geometrical realization. It is then natural to ask to what extend
   the combinatorial upper bound on the number of triple points
   found by Sch\"onhein \cite{Sch66} (see also equation \eqnref{eq: schoenheim} below)
   can be improved in general.
   Our result is rendered by a series of configurations
   with a high number of triple points presented in Section \ref{sec:triple configs}.
   This series of examples sets some limits to possible improvements in \eqnref{eq: schoenheim}.

   Given a positive integer $s$ and a projective plane over a field with
   sufficiently many elements, it is easy to find $s$ lines intersecting in
   exactly $\binom{s}{2}$ distinct points. In fact this is the number of
   intersection points of a \emph{general} arrangement of
   $s$ lines (such arrangements in algebraic geometry are called
   \emph{star configurations}, see \cite{GHM13}) and this is also the maximal possible number of points
   in which at least two out of given $s$ lines intersect.

   Given a configuration of $s$ mutually distinct lines, let $t_k$
   denote the number of points where exactly $k\geq 2$ lines meet.
   Then there is the obvious combinatorial equality
   \begin{equation}\label{eq: combinatorial}
      \binom{s}{2}=\sum\limits_{k\geq 2}t_k\binom{k}{2}.
   \end{equation}
   Let $T_k(s)$ denote the maximal number of $k$-fold points
   in an arrangement of $s$ distinct lines in the projective plane
   over an arbitrary field $\F$. The discussion above shows that
   $$T_2(s)=\binom{s}{2}.$$
   The aim of this note is to investigate the numbers $T_3(s)$.

\section{Arrangements with many triple points}
   Deciding the existence or non-existence of a configuration with
   certain properties is a problem which can, in principle, be always
   solved by combinatorial and semi-algebraic methods. The combinatorial
   part evaluates the collinearity conditions and checks whether
   the resulting incidence table can be filled in or not. The restrictions
   are imposed by the number of lines intersecting in configuration points
   and the condition that two lines cannot intersect in more than one point.
   Table \ref{tab: inc 10 lines} on page \pageref{tab: M1234} is an example of what
   we call an incidence table.

   If a configuration is combinatorially possible, then we assign
   coordinates to the equations of configuration lines and check if the
   system of polynomial equations resulting from evaluating
   combinatorial data has solutions. Typically this is the case
   and this is where the semi-algebraic part comes into the picture,
   as one has to exclude various degenerations, for example
   points or lines falling together. This kind of conditions is
   given by \emph{inequalities}.

   This note, in a sense, is a field case study of the effective
   applicability of the approach described above. Evaluating
   all conditions carefully, one can actually study \emph{moduli spaces}
   of configurations in the spirit of \cite{ATY13} and \cite{ACTY}.
   However, since we are interested in the existence of configurations
   over \emph{arbitrary} fields, we do not dwell on this aspect
   of the story.

   Some of the computations
   were supported by the symbolic algebra program Singular \cite{DGPS}.

   The equality \eqnref{eq: combinatorial} yields the following
   upper bound on the number $t_3$ of triple points in an arrangement
   of $s$ lines:
   \begin{equation}\label{eq: naive}
      t_3\leq \rounddown{\frac{\binom{s}{2}}{3}}.
   \end{equation}
   This naive bound has been improved by
   Kirkman in 1847, \cite{Kir47}, with a correction
   of Sch\"onheim in 1966, \cite{Sch66}. Theorem \ref{thm: main2} shows that
   this bound is close to be attained, on the other hand Theorem \ref{thm: main}
   shows that there is place for some further improvements. Let
   \begin{equation}\label{eq: schoenheim}
      U_3(s):=\rounddown{\rounddown{\frac{s-1}{2}}\cdot\frac{s}{3}}-\eps(s),
   \end{equation}
   where $\eps(s)=1$ if $s \equiv 5 \mod(6)$ and $\eps(s)=0$ otherwise. Then
   \begin{equation}\label{eq: upper on 3ple points}
      T_3(s)\leq U_3(s).
   \end{equation}
   We refer to Section 5 in \cite{BGS74} for a nice discussion of historical
   backgrounds. In the next table we present a few first numbers resulting from
   \eqnref{eq: upper on 3ple points}
   \begin{center}
   \renewcommand{\arraystretch}{1.2}
   \begin{tabular}{|c|c|c|c|c|c|c|c|c|c|c|c|c|}
   \hline
   $s$ & 1 & 2 & 3 & 4 & 5 & 6 & 7 & 8 & 9 & 10 & 11 & 12 \\
   \hline
   $U_3(s)$ & 0 & 0 & 1 & 1 & 2 & 4 & 7 & 8 & 12 & 13 & 17 & 20  \\
   \hline
   \end{tabular}
   \end{center}
   It is natural to ask to which extend the numbers appearing in the above
   table are sharp.
   Our main result is the following classification Theorem.
\begin{theorem}\label{thm: main}
\begin{itemize}
   \item[a)] For $1\leq s\leq 6$, there are configurations of lines with $U_3(s)$ triple points
   in projective planes $\P(\F)$ over arbitrary fields.
   \item[b)] A configuration of $7$ lines with $7$ triple points exists only in characteristic $2$
   (the smallest such configuration is the Fano plane $\P^2(\F_2)$).
   \item[c)] A configuration of $8$ lines with $8$ triple points exists over any field
   containing a non-trivial third degree root of $1$. Moreover such a configuration
   always arises by taking out one line from a configuration of $9$ lines
   with $12$ triple points.
   \item[d)] A configuration of $9$ lines with $12$ triple points exists over
   any field containing a non-trivial third degree root of $1$.
   \item[e)] A configuration of $10$ lines with $13$ triple points exists
   only
   \subitem e1) over a field $\F$ of characteristic $2$ containing a non-trivial third root of unity. In this case
   one of the points has in fact multiplicity $4$;
   \subitem e2) over any field $\F$ of characteristic $5$.
   \item[f)] There is no configuration of $11$ lines with $17$ triple points.
   There exist configurations of $11$ lines with $16$ triple points.
   \end{itemize}
\end{theorem}
\proof
   The first part of the Theorem is well known for $s\leq 9$.
   We go briefly through
   all the cases for the sake of the completeness and discuss
   $s=10$ in more detail as this configuration
   seems to be new.

   For $s=1,2$ there
   is nothing to prove. For $s=3,4$ we take $3$ lines in a pencil
   and an arbitrary fourth line. The case $s=5$ is easy as well, see Figure \ref{fig: s=5}.
\begin{figure}[H]
\centering
   \begin{minipage}{0.4\textwidth}
   \centering
\begin{tikzpicture}[line cap=round,line join=round,x=1.0cm,y=1.0cm,scale=0.7]
\clip(1.62,-2.46) rectangle (7.72,2.74);
\draw [domain=1.62:7.72] plot(\x,{(-12.49--4.64*\x)/5.34});
\draw [domain=1.62:7.72] plot(\x,{(-6.19--2.2*\x)/-1.92});
\draw [domain=1.62:7.72] plot(\x,{(--0.17--0.02*\x)/3.82});
\draw [domain=1.62:7.72] plot(\x,{(-27.01--4.06*\x)/-3.74});
\draw [domain=1.62:7.72] plot(\x,{(-29.85--4.6*\x)/5.24});
\begin{scriptsize}
\fill [color=black] (-2.58,-4.58) circle (2.5pt);
\draw[color=black] (-2.44,-4.3) node {$A$};
\fill [color=black] (2.76,0.06) circle (2.5pt);
\fill [color=black] (0.84,2.26) circle (2.5pt);
\draw[color=black] (1,2.54) node {$C$};
\fill [color=black] (6.58,0.08) circle (2.5pt);
\fill [color=black] (2.84,4.14) circle (2.5pt);
\draw[color=black] (3,4.42) node {$E$};
\fill [color=black] (11.82,4.68) circle (2.5pt);
\draw[color=black] (11.96,4.96) node {$F$};
\end{scriptsize}
\end{tikzpicture}
   \caption{$ $ : $s=5$}\label{fig: s=5}
   \end{minipage}
   \quad
   \begin{minipage}{0.4\textwidth}
   \centering
\begin{tikzpicture}[line cap=round,line join=round,x=1.0cm,y=1.0cm,scale=0.7]
\clip(1.62,-2.46) rectangle (7.72,2.74);
\draw [domain=1.96:7.74] plot(\x,{(-12.49--4.64*\x)/5.34});
\draw [domain=1.96:7.74] plot(\x,{(-6.19--2.2*\x)/-1.92});
\draw [domain=1.96:7.74] plot(\x,{(--0.17--0.02*\x)/3.82});
\draw [domain=1.96:7.74] plot(\x,{(-27.01--4.06*\x)/-3.74});
\draw [domain=1.96:7.74] plot(\x,{(-29.85--4.6*\x)/5.24});
\draw [domain=1.96:7.74] plot(\x,{(--17.37-3.74*\x)/-0.48});
\begin{scriptsize}
\fill [color=black] (-2.58,-4.58) circle (1.5pt);
\draw[color=black] (-2.44,-4.3) node {$A$};
\fill [color=black] (2.76,0.06) circle (2.5pt);
\fill [color=black] (0.84,2.26) circle (1.5pt);
\draw[color=black] (1,2.54) node {$C$};
\fill [color=black] (6.58,0.08) circle (2.5pt);
\fill [color=black] (2.84,4.14) circle (1.5pt);
\draw[color=black] (3,4.42) node {$E$};
\fill [color=black] (11.82,4.68) circle (1.5pt);
\draw[color=black] (11.96,4.96) node {$F$};
\fill [color=black] (4.89,1.91) circle (2.5pt);
\fill [color=black] (4.41,-1.83) circle (2.5pt);
\end{scriptsize}
\end{tikzpicture}
\caption{$ $ : $s=6$}\label{fig: s=6}
\end{minipage}
\end{figure}
   We pass to the case $s=6$ adding the line through both
   double points, see Figure~\ref{fig: s=6}.

   For $s=7$ we obtain the famous Fano plane
   $\P^2(\F_2)$. It is well known that this configuration is possible
   only in characteristic $2$. The picture below (Figure \ref{fig: s=7})
   indicates collinear points as lying on the
   segments or on the circle.
\begin{figure}[H]
\centering
\begin{minipage}{0.4\textwidth}
\centering
\begin{tikzpicture}[line cap=round,line join=round,x=1.0cm,y=1.0cm,scale=0.7]
\clip(0.02,-2.74) rectangle (8,3.34);
\draw [color=zzttqq] (1.74,-1.7)-- (6.72,-1.7);
\draw [color=zzttqq] (6.72,-1.7)-- (4.23,2.61);
\draw
 [color=zzttqq] (4.23,2.61)-- (1.74,-1.7);
\draw(4.23,-0.26) circle (1.44cm);
\draw (1.74,-1.7)-- (5.48,0.46);
\draw (4.23,2.61)-- (4.23,-1.7);
\draw (2.99,0.46)-- (6.72,-1.7);
\begin{scriptsize}
\fill [color=qqqqff] (1.74,-1.7) circle (2.5pt);
\fill [color=qqqqff] (6.72,-1.7) circle (2.5pt);
\fill [color=uuuuuu] (4.23,2.61) circle (2.5pt);
\fill [color=uuuuuu] (4.23,-1.7) circle (2.5pt);
\fill [color=uuuuuu] (2.99,0.46) circle (2.5pt);
\fill [color=uuuuuu] (5.48,0.46) circle (2.5pt);
\fill [color=uuuuuu] (4.23,-0.26) circle (2.5pt);
\end{scriptsize}
\end{tikzpicture}
\caption{$ $ : $s=7$}
\label{fig: s=7}
\end{minipage}
\quad
\begin{minipage}{0.4\textwidth}
   \centering
\begin{tikzpicture}[line cap=round,line join=round,x=1.0cm,y=1.0cm,scale=0.7]
\clip(0.02,-2.74) rectangle (8,3.34);
\draw [color=zzttqq] (1.74,-1.7)-- (6.72,-1.7);
\draw [color=zzttqq] (6.72,-1.7)-- (4.23,2.61);
\draw [color=zzttqq] (4.23,2.61)-- (1.74,-1.7);
\draw (1.74,-1.7)-- (5.48,0.46);
\draw (2.99,0.46)-- (6.72,-1.7);
\draw (2.99,0.46)-- (4.23,-1.7);
\draw (4.23,-1.7)-- (5.48,0.46);
\draw [shift={(4.23,0.94)}] plot[domain=-1.19:4.33,variable=\t]({1*1.68*cos(\t r)+0*1.68*sin(\t r)},{0*1.68*cos(\t r)+1*1.68*sin(\t r)});
\begin{scriptsize}
\fill [color=qqqqff] (1.74,-1.7) circle (2.5pt);
\fill [color=qqqqff] (6.72,-1.7) circle (2.5pt);
\fill [color=uuuuuu] (4.23,2.61) circle (2.5pt);
\fill [color=uuuuuu] (4.23,-1.7) circle (2.5pt);
\fill [color=uuuuuu] (2.99,0.46) circle (2.5pt);
\fill [color=uuuuuu] (5.48,0.46) circle (2.5pt);
\fill [color=uuuuuu] (3.61,-0.62) circle (2.5pt);
\fill [color=uuuuuu] (4.85,-0.62) circle (2.5pt);
\end{scriptsize}
\end{tikzpicture}
\caption{$ $ : $s=8$}
\label{fig: s=8}
\end{minipage}
\end{figure}
   For $s=8$ there is the M\"obius-Kantor $(8_3)$ configuration. This
   configuration cannot be drawn in the real plane. Collinearity
   is indicated by segments and the circle arch, see Figure \ref{fig: s=8}. This configuration
   can be obtained from the next configuration by removing one line.

   For $s=9$ there is the dual Hesse configuration.
   It is easier to describe the original
   Hesse configuration. It arises taking the nine
   order $3$ torsion points of a smooth
   complex cubic curve (which carries the structure
   of an abelian group and the torsion is understood
   with respect to this group structure). There are
   $12$ lines passing through the nine points in such
   a way that each line contains exactly $3$ torsion
   points and there are $4$ lines passing through each
   of the points. See \cite{ArtDol09} for details.
   This configuration cannot
   be drawn in the real plane.\\
   Beside the geometrical realization over the complex numbers,
   the dual Hesse configuration can be also easily obtained
   in characteristic $3$, more precisely,
   in the plane $\P^2(\F_3)$ taking all $13$ lines and removing
   from this set all $4$ lines passing through a fixed point.
   We leave the details to the reader.

   Before moving on, we record for further reference the following
   simple but useful fact.
\begin{lemma}\label{lem: points on one line}
   Let $\call=\left\{L_1,\ldots,L_s\right\}$ be a configuration of lines.
   Let $L\in\call$ be a fixed line. Let $P_1(L),\ldots,P_r(L)$
   be intersection points of $L$ with other configuration lines
   with corresponding multiplicities $m_1(L),\ldots,m_r(L)$. Then
   $$s-1=\sum\limits_{i=1}^r(m_i-1).$$
   In particular, if there are only triple points on a line $L$,
   then $s$ is an odd number.
\end{lemma}
\subsection{$10$ lines}
   Now we come to the case $s=10$. We work over an arbitrary field $\F$.
   The upper bound \eqnref{eq: upper on 3ple points} implies that in this case
   there can be at most $13$ points of multiplicity at least $3$. We first
   deal with the case when points of higher multiplicity might appear.
\subsubsection{Points of excess multiplicity}
   The combinatorial equality \eqnref{eq: combinatorial}
   implies that there are no points with multiplicity $m\geq 5$ and
   only the following cases with $4$-fold points need to be considered:
\begin{enumerate}
\item[(i)] $t_4=2$, $t_3=11$, $t_2=0$.
\item[(ii)] $t_4=1$, $t_3=12$, $t_2=3$.
\end{enumerate}
   Case (i) is excluded since there exists a configuration line $L$ which does not pass through any of the $4$-fold points. Hence this line contains only $3$-fold points. Since $s=10$, this contradicts Lemma \ref{lem: points on one line}.

   Passing to (ii) we start with the $4$-fold point $W$. We denote the lines passing through $W$
   by $M_1,\ldots,M_4$ as indicated in the picture below.
\begin{figure}[H]
\centering
\begin{tikzpicture}[line cap=round,line join=round,x=0.81cm,y=0.81cm, scale=0.7]
\clip(-7.0,-2.06) rectangle (6.0,8.39);
\draw [domain=-7.0:4.0] plot(\x,{(--21.03--3.296*\x)/7.39}); 
\draw [domain=-7.0:1.0] plot(\x,{(--20.86--4.51*\x)/1.472});
\draw [domain=-7.0:1.5] plot(\x,{(--23.46--4.46*\x)/4.53}); 
\draw [domain=-7.0:5.0] plot(\x,{(--8.34--0.66*\x)/6.026}); 
\begin{scriptsize}
\draw [fill=black] (-4.33,0.914) circle (2.5pt);
\draw[color=black] (-4.486,1.24) node {$W$};
\draw [fill=black] (3.06,4.21) circle (2.5pt);
\draw[color=black] (2.994,4.672) node {$P_{35}$};
\draw[color=black] (-6.488,0.323) node {$M_3$};
\draw [fill=black] (-2.858,5.42) circle (2.5pt);
\draw[color=black] (-3.21,5.86) node {$P_{34}$};
\draw[color=black] (-5.6,-1.516) node {$M_1$};
\draw [fill=black] (0.2,5.376) circle (2.5pt);
\draw[color=black] (-0.086,5.86) node {$P_{46}$};
\draw[color=black] (-6.578,-0.698) node {$M_2$};
\draw [fill=black] (1.696,1.57) circle (2.5pt);
\draw[color=black] (1.674,1.944) node {$P_{24}$};
\draw[color=black] (-6.488,1.042) node {$M_4$};
\draw [fill=black] (-3.2533470849750925,4.2097867086292355) circle (2.5pt);
\draw[color=black] (-3.584,4.694) node {$P_{12}$};
\draw [fill=black] (-2.4369368932686926,6.708933667752222) circle (2.5pt);
\draw[color=black] (-2.748,7.2) node {$P_{56}$};
\draw [fill=black] (-2.143711375856662,3.0674701635601713) circle (2.5pt);
\draw[color=black] (-2.352,3.55) node {$P_{13}$};
\draw [fill=black] (-0.8614830881576307,4.33045087431361) circle (2.5pt);
\draw[color=black] (-1.142,4.782) node {$P_{25}$};
\draw [fill=black] (-0.7212666479632315,2.5235243746025966) circle (2.5pt);
\draw[color=black] (-0.812,3.0) node {$P_{14}$};
\draw [fill=black] (1.2038034131825759,3.3821212516711463) circle (2.5pt);
\draw[color=black] (1.124,3.858) node {$P_{26}$};
\draw [fill=black] (-0.471031159814296,1.3340935212681417) circle (2.5pt);
\draw[color=black] (-0.46,1.68) node {$P_{15}$};
\draw [fill=black] (3.9795739846919567,1.8185935170856165) circle (2.5pt);
\draw[color=black] (3.984,2.186) node {$P_{36}$};
\end{scriptsize}
\end{tikzpicture}
\caption{}
\label{fig: with W}
\end{figure}

   The remaining $6$ lines $L_1,\ldots, L_6$ (not visible in the figure above) intersect pairwise in $15$ mutually distinct points --- $12$ of these points are the $12$ configuration triple points and the remaining $3$ points (not visible in the figure above) contribute to $t_2$.
   Note that $6$ general lines intersect in $15$ mutually distinct points,
   but in our situation there are additional collinearities which are reflected in the picture above and in the table
   below. Passing to the details let $P_{ij}=L_i\cap L_j$ for $1\leq i<j\leq 6$.
   Up to renumbering of points we may assume that they are distributed in the following way:
   \begin{table}[H]
   \centering
   \renewcommand{\arraystretch}{1.2}
   \begin{tabular}{|c|c|}
   \hline
   line & points on the line \\
   \hline
   $M_1$ & $P_{12}$, $P_{34}$, $P_{56}$  \\
   \hline
   $M_2$ & $P_{13}$, $P_{25}$, $P_{46}$  \\
   \hline
   $M_3$ & $P_{14}$, $P_{26}$, $P_{35}$  \\
   \hline
   $M_4$ & $P_{15}$, $P_{24}$, $P_{36}$  \\
   \hline
   \end{tabular}
   \caption{}
   \label{tab: M1234}
   \end{table}

   Moreover, we may assume that the lines $L_1,\ldots, L_6$ have the following equations (we omit ``$=0$''   in the equations)
   \begin{equation}\label{eq: abcd}
   \begin{array}{lll}
   L_1:\; x,          & L_2:\; y,        & L_3:\; z,     \\
   L_4:\; x+y+z,      & L_5:\; ax+by+z,  & L_6:\; cx+dy+z,
   \end{array}
   \end{equation}
   with
   \begin{equation}\label{eq: abcd det}
   a,b,c,d\in\mathbb{F}^*\;\;\mbox{ and }\; \det\left(\begin{matrix}1&1&1\\a&b&1\\c&d&1\end{matrix}\right)\neq 0.
   \end{equation}
   Indeed, the first four equations are obvious. The coefficients at $z$ in the lines $L_5$ and $L_6$ can be
   normalized to $1$ since otherwise the star configuration condition would fail.
   Similarly the conditions in \eqnref{eq: abcd det} are necessary in order to
   guarantee that $L_1,\ldots,L_6$ form a star configuration.

   Evaluating collinearity conditions in Table \ref{tab: M1234} above
   we obtain the following system of linear and quadratic equations:
\begin{equation}\label{eq: 10-1-12-3_1}
\left\{\begin{array}{r} a-b-c+d=0\\-ad+a-c+d=0\\a-bc=0\\bc-d=0\end{array}\right.
\end{equation}

   Additionally, the condition that the lines $M_1, \ldots, M_4$ belong to the same pencil gives
\begin{equation}\label{eq: 10-1-12-3_2}
-ab+a+bc-1=0.
\end{equation}

   A solution to the above system of equations
   \eqnref{eq: 10-1-12-3_1} and \eqnref{eq: 10-1-12-3_2}
   satisfying additionally the non-equality condition \eqnref{eq: abcd det}
   exists only in characteristic $2$. This has been verified with the aid of Singular.
   Moreover, in that case
   $a$ satisfies
   $$a^2+a+1=0.$$
   and then, consequently, $b=a^2$, $c=a^2$ and $d=a$.

   It follows that the configuration (ii) exists in $\P^2(\mathbb{F}_{2^q})$, for all $q\geq 2$
   (the case $q=1$ is excluded
   as there are evidently not enough points in the Fano plane.
   The configuration lines are then given by equations
   \begin{equation*}
   \begin{array}{llll}
    L_1:\; x,          & L_2:\; y,        & L_3:\; z,     &\\
    L_4:\; x+y+z,     & L_5:\; ax+a^2y+z, & L_6:\; a^2x+ay+z,   & \\
    M_1:\;x+y ,   & M_2:\; ax+z,    & M_3:\; a^2x+y+z,      & M_4:\; x+a^2y+z
   \end{array}
   \end{equation*}

   Then the configuration points have coordinates:
  \begin{equation*}
   \begin{array}{llll}
      W=(1:1:a),              & P_{12}=(0:0:1),         & P_{13}=(0:1:0),        & P_{14}=(0:1:1),\\
      P_{15}=(0:1:a^2),       & P_{24}=(1:0:1),         & P_{25}=(1:0:a),        & P_{26}=(1:0:a^2),\\
      P_{34}=(1:1:0),         & P_{35}=(a:1:0),         & P_{36}=(1:a:0),        & P_{46}=(a^2:a:1),\\
      P_{56}=(1:1:1).\\
   \end{array}
   \end{equation*}

   The incidence table in this case reads:\label{page:tab}
   \begin{table}[H]
   \centering
   \renewcommand{\arraystretch}{1.2}
   \begin{tabular}{|c|c|c|c|c|c|c|c|c|c|c|c|c|c|}
   \hline
   & $W$ & $P_{12}$ & $P_{13}$ & $P_{14}$ & $P_{15}$ & $P_{24}$ & $P_{25}$ & $P_{26}$ & $P_{34}$ & $P_{35}$ & $P_{36}$ & $P_{46}$ & $P_{56}$ \\
   \hline
   $L_1$ &  &+ &+ &+ &+ &  &  &  &  &  &  &  & \\
   \hline
   $L_2$ &  &+ &  &  &  &+ &+ &+ &  &  &  &  & \\
   \hline
   $L_3$ &  &  &+ &  &  &  &  &  &+ &+ &+ &  & \\
   \hline
   $L_4$ &  &  &  &+ &  &+ &  &  &+ &  &  &+ & \\
   \hline
   $L_5$ &  &  &  &  &+ &  &+ &  &  &+ &  &  &+\\
   \hline
   $L_6$ &  &  &  &  &  &  &  &+ &  &  &+ &+ &+\\
   \hline
   $M_1$ &+ &+ &  &  &  &  &  &  &+ &  &  &  &+\\
   \hline
   $M_2$ &+ &  &+ &  &  &  &+ &  &  &  &  &+ & \\
   \hline
   $M_3$ &+ &  &  &+ &  &  &  &+ &  &+ &  &  & \\
   \hline
   $M_4$ &+ &  &  &  &+ &+ &  &  &  &  &+ &  & \\
   \hline
   \end{tabular}
   \caption{}
   \label{tab: inc 10 lines}
   \end{table}
\subsubsection{Points of multiplicity $3$}
   Now we pass to the case that there are no quadruple points, hence $t_3=13$ and consequently $t_2=6$.
   There is an odd number of $2$-fold points on each configuration line. This implies that there is a configuration line $M_1$ containing exactly three $2$-fold points $D_1$, $D_2$, $D_3$. We have again two cases:
\begin{enumerate}
\item[(A)] The lines $M_2$, $M_3$, $M_4$ passing through the points $D_1$, $D_2$, $D_3$ meet in a single point $W$.
\begin{figure}[H]
\centering
\begin{tikzpicture}[line cap=round,line join=round,x=1.0cm,y=1.0cm, scale=0.7]
\clip(-4.3,0.5) rectangle (4.5,6.3);
\draw [domain=-4.3:4.5] plot(\x,{(--12.8152-0.019999999999999574*\x)/2.88});
\draw [domain=-4.3:4.5] plot(\x,{(--3.0204-2.5999999999999996*\x)/1.54});
\draw [domain=-4.3:4.5] plot(\x,{(--0.3789941071514682-2.5901403288807447*\x)/0.12020735882721695});
\draw [domain=-4.3:4.5] plot(\x,{(-2.3376-2.58*\x)/-1.3399999999999999});
\begin{scriptsize}
\draw [fill=qqqqff] (-1.48,4.46) circle (2.5pt);
\draw[color=qqqqff] (-1.25,4.8) node {$D_1$};
\draw [fill=qqqqff] (1.4,4.44) circle (2.5pt);
\draw[color=qqqqff] (1.15,4.8) node {$D_3$};
\draw[color=black] (-3.55,4.15) node {$M_1$};
\draw [fill=qqqqff] (-0.06020735882721695,4.450140328880745) circle (2.5pt);
\draw[color=qqqqff] (0.27,4.8) node {$D_2$};
\draw [fill=qqqqff] (0.06,1.86) circle (2.5pt);
\draw[color=qqqqff] (0.55,1.95) node {$W$};
\draw[color=black] (-1.95,6.0) node {$M_2$};
\draw[color=black] (0.3,6.0) node {$M_3$};
\draw[color=black] (1.7,6.0) node {$M_4$};
\end{scriptsize}
\end{tikzpicture}
\caption{$ $ : Case (A)}
\label{fig: case A}
\end{figure}
\item[(B)] The lines $M_2$, $M_3$, $M_4$ form a triangle with vertices $Z_1$, $Z_2$, $Z_3$.
\begin{figure}[H]
\centering
\begin{tikzpicture}[line cap=round,line join=round,x=1.0cm,y=1.0cm,scale=0.7]
\clip(-3.3,0.7) rectangle (5.4,6.3);
\draw [domain=-3.3:5.4] plot(\x,{(--14.21-0.0*\x)/3.32});
\draw [domain=-3.3:5.4] plot(\x,{(--12.01-1.44*\x)/2.84});
\draw [domain=-3.3:5.4] plot(\x,{(--4.62-2.66*\x)/0.06});
\draw [domain=-3.3:5.4] plot(\x,{(--2.06-2.66*\x)/-1.52});
\begin{scriptsize}
\draw [fill=qqqqff] (-0.1,4.28) circle (2.5pt);
\draw [fill=qqqqff] (3.22,4.28) circle (2.5pt);
\draw[color=black] (-2.8,3.9) node {$M_1$};
\draw [fill=qqqqff] (1.64,4.28) circle (2.5pt);
\draw[color=black] (-2.8,5.3) node {$M_2$};
\draw[color=black] (1.15,6.0) node {$M_3$};
\draw[color=black] (3.64,6.0) node {$M_4$};
\draw [fill=uuuuuu] (1.66,3.39) circle (2.5pt);
\draw [fill=uuuuuu] (2.47,2.97) circle (2.5pt);
\draw [fill=uuuuuu] (1.70,1.62) circle (2.5pt);
\draw[color=black] (1.2,1.78) node {$Z_1$};
\draw[color=black] (3.0,3.0) node {$Z_2$};
\draw[color=black] (1.18,3.2) node {$Z_3$};
\end{scriptsize}
\end{tikzpicture}
\caption{$ $ : Case (B)}
\label{fig: case B}
\end{figure}
\end{enumerate}

The first case is impossible. This follows similarly to the case (e1) with a $4$-fold point.
Indeed, the six lines $L_1, \ldots, L_6$ not visible in the Figure \ref{fig: case A}
form a star configuration, i.e. there are 15 mutually
distinct intersection points $P_{ij}=L_i\cap L_j$. Up to renumbering incidences between the lines
$M_i$ and the points $P_{jk}$ are as in the Table \ref{tab: M1234 2}.

Moreover the condition that $M_2$, $M_3$ and $M_4$ meet at one point $W$ is the same as in equation \eqnref{eq: 10-1-12-3_2}.
This gives the same solution as in the previous case (e1). It is easy to check that this implies that the line $M_1$ goes through $W$, a contradiction.

   Now we consider the remaining case (B). There are three $3$-fold points on the line $M_1$. The
points $Z_1$, $Z_2$, $Z_3$ must be also $3$-fold points of the configuration and on each of the
lines $M_2$, $M_3$, $M_4$ there are two more $3$-fold points. This gives altogether twelve $3$-fold
points. Hence the six remaining lines $L_1, \ldots, L_6$ have a $3$-fold intersection point. We
call this point $D$, and assume that $L_4$, $L_5$ and $L_6$ pass through $D$. We can assume that
   $D=(1:1:1)$ and then the
equations of the lines $L_i$ are
\begin{equation}\label{eq: 10-13_3}
   \begin{array}{lll}
   L_1:\; x,          & L_2:\; y,        & L_3:\; z,     \\
   L_4:\; ax-(a+1)y+z,     & L_5:\; bx-(b+1)y+z, & L_6:\; cx-(c+1)y+z,
   \end{array}
   \end{equation}
with $a$, $b$, $c$ mutually distinct and different from zero. The combinatorics implies that
each of $Z_1$, $Z_2$, $Z_3$ lie on one of the lines
$L_4$, $L_5$, $L_6$. Up to renumbering we can assume $Z_1\in L_4$, $Z_2\in L_5$, $Z_3\in L_6$. Then
the incidence table is determined as follows:

\begin{table}[H]
   \centering
   \renewcommand{\arraystretch}{1.2}
   \begin{tabular}{|c|c|}
   \hline
   line & points on the line \\
   \hline
   $M_1$ & $P_{14}$, $P_{25}$, $P_{36}$  \\
   \hline
   $M_2$ & $Z_{2}$, $Z_{3}$, $P_{12}$, $P_{34}$  \\
   \hline
   $M_3$ & $Z_{1}$, $Z_{3}$, $P_{15}$, $P_{23}$  \\
   \hline
   $M_4$ & $Z_{1}$, $Z_{2}$, $P_{13}$, $P_{26}$ \\
   \hline
   \end{tabular}
   \caption{}
   \label{tab: M1234 2}
\end{table}
   See the text after Table \ref{tab: NiMi} for hints how to fill in such a table.

Using equations as in equation \eqnref{eq: 10-13_3} and evaluating incidences we obtain the following conditions
\begin{equation}
\left\{\begin{array}{r} ab+ac+a-bc=0\\ac+a-b+c=0\\ab+c=0\\a+bc=0\end{array}\right.
\end{equation}

This implies that $b^2=1$. If $b=-1$ then $a=-1$, a contradiction. If $b=1$ then $a=3$ and $a^2+1=0$. It follows that $\characteristic \mathbb{F}=2$ or $\characteristic \mathbb{F}=5$. In the first case $a=c$, a contradiction. In the second case we obtain $a=3$, $b=1$, $c=2$, thus the lines are
   \begin{equation*}
   \begin{array}{llll}
    L_1:\; x,          & L_2:\; y,           & L_3:\; z, &\\
    L_4:\; 3x+y+z,     & L_5:\; x+3y+z, &    L_6:\; 2x+2y+z,      & \\
    M_1:\;x+y+z ,   & M_2:\; 2x+4y,    & M_3:\; 3y+z,      & M_4:\; 2x+z
   \end{array}
   \end{equation*}

The points have coordinates
  \begin{equation*}
   \begin{array}{llll}
      D=(1:1:1),              & Z_1=(2:3:1),          & Z_2=(4:3:2),         & Z_3=(4:3:1),\\
      P_{12}=(0:0:1),         & P_{13}=(0:1:0),        & P_{14}=(0:4:1),     & P_{15}=(0:4:3),\\
      P_{23}=(1:0:0),         & P_{25}=(1:0:4),        & P_{26}=(1:0:3),      & P_{34}=(4:3:0),\\
      P_{36}=(3:2:0).
   \end{array}
   \end{equation*}

The incidence table is
\begin{table}[H]
   \centering
   \renewcommand{\arraystretch}{1.2}
   \begin{tabular}{|c|c|c|c|c|c|c|c|c|c|c|c|c|c|}
   \hline
   & $D$& $Z_1$ & $Z_2$  & $Z_3$ & $P_{12}$ & $P_{13}$ & $P_{14}$ & $P_{15}$ & $P_{23}$ & $P_{25}$ & $P_{26}$ & $P_{34}$ & $P_{36}$ \\
   \hline
   $L_1$ &  &  &  &  &+ &+ &+ &+ &  &  &  &  & \\
   \hline
   $L_2$ &  &  &  &  &+ &  &  &  &+ &+ &+ &  & \\
   \hline
   $L_3$ &  &  &  &  &  &+ &  &  &+ &  &  &+ &+\\
   \hline
   $L_4$ &+ &+ &  &  &  &  &+ &  &  &  &  &+ & \\
   \hline
   $L_5$ &+ &  &+ &  &  &  &  &+ &  &+ &  &  & \\
   \hline
   $L_6$ &+ &  &  &+ &  &  &  &  &  &  &+ &  &+\\
   \hline
   $M_1$ &  &  &  &  &  &  &+ &  &  &+ &  &  &+\\
   \hline
   $M_2$ &  &  &+ &+ &+ &  &  &  &  &  &  &+ & \\
   \hline
   $M_3$ &  &+ &  &+ &  &  &  &+ &  &  &+ &  & \\
   \hline
   $M_4$ &  &+ &+ &  &  &+ &  &  &+ &  &  &  & \\
   \hline
   \end{tabular}
   \caption{}
   \label{tab: inc 10 lines 2}
\end{table}

Passing to the last assertion
 f) of Theorem \ref{thm: main} we assume that a configuration of 11 lines with 17 triple points exists. Then \eqnref{eq: combinatorial} implies that there are 4 double points in the configuration. Hence each line meets 10 other lines, there is an even number of double points on each line. If there are 4 double points on a line, then this condition fails on the 4 lines meeting the given one in the double points. Hence, there must be 2 pairs of lines in the configuration with double points situated in the intersection points of lines from different pairs as indicated in the figure below
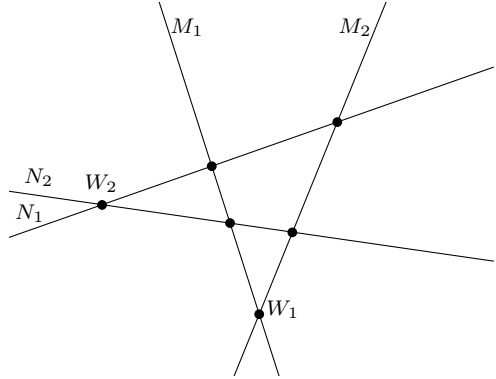
\begin{center}
\begin{tikzpicture}[line cap=round,line join=round,x=1.0cm,y=1.0cm, scale=0.65]
\clip(-4.3,-1.3) rectangle (5.5,6.3);
\draw [domain=-4.3:5.5] plot(\x,{(--7.7--0.9*\x)/2.54});
\draw [domain=-4.3:5.5] plot(\x,{(--2.24-3.02*\x)/0.96});
\draw [domain=-4.3:7.06] plot(\x,{(--3.074-3.92*\x)/-1.58});
\draw [domain=-4.3:5.5] plot(\x,{(--7.02-0.56*\x)/3.85});
\begin{scriptsize}
\draw [fill=qqqqff] (-0.2,2.96) circle (2.5pt);
\draw [fill=qqqqff] (2.34,3.86) circle (2.5pt);
\draw[color=black] (-3.88,2.0) node {$N_1$};
\draw [fill=qqqqff] (0.76,-0.06) circle (2.5pt);
\draw[color=qqqqff] (1.24,0.04) node {$W_1$};
\draw[color=black] (-0.7,5.8) node {$M_1$};
\draw[color=black] (2.7,5.8) node {$M_2$};
\draw [fill=qqqqff] (-2.42,2.17) circle (2.5pt);
\draw[color=qqqqff] (-2.44,2.62) node {$W_2$};
\draw [fill=qqqqff] (1.43,1.61) circle (2.5pt);
\draw[color=black] (-3.7,2.74) node {$N_2$};
\draw [fill=qqqqff] (0.17,1.8) circle (2.5pt);
\end{scriptsize}
\end{tikzpicture}
\end{center}

Now, there are two cases
\begin{enumerate}
\item[(I)] the line $W_1W_2$ belongs to the configuration,
\item[(II)]  the line $W_1W_2$ is not a  configuration line.
\end{enumerate}

We begin with the case (I)
\begin{figure}[H]
\centering
\begin{tikzpicture}[line cap=round,line join=round,x=1.0cm,y=1.0cm, scale=0.6]
\clip(-5,-2.0) rectangle (5.5,6.3);
\draw [domain=-5:5.5] plot(\x,{(--7.7--0.9*\x)/2.54});
\draw [domain=-5:5.5] plot(\x,{(--2.24-3.02*\x)/0.96});
\draw [domain=-5:5.5] plot(\x,{(--3.074-3.92*\x)/-1.58});
\draw [domain=-5:5.5] plot(\x,{(--7.02-0.56*\x)/3.85});
\draw [domain=-5:5.5] plot(\x,{(--1.51-2.23*\x)/3.18});
\begin{scriptsize}
\draw [fill=qqqqff] (-0.2,2.96) circle (2.5pt);
\draw [fill=qqqqff] (2.34,3.86) circle (2.5pt);
\draw[color=black] (-4.3,1.85) node {$N_1$};
\draw [fill=qqqqff] (0.76,-0.06) circle (2.5pt);
\draw[color=qqqqff] (1.4,0.0) node {$W_1$};
\draw[color=black] (-0.62,5.8) node {$M_1$};
\draw[color=black] (2.65,5.8) node {$M_2$};
\draw [fill=qqqqff] (-2.42,2.17) circle (2.5pt);
\draw[color=qqqqff] (-2.4,2.62) node {$W_2$};
\draw [fill=qqqqff] (1.43,1.61) circle (2.5pt);
\draw[color=black] (-4.3,2.74) node {$N_2$};
\draw[color=black] (-4.43,3.95) node {$L$};
\draw [fill=qqqqff] (0.17,1.8) circle (2.5pt);
\end{scriptsize}
\end{tikzpicture}
\caption{$ $ : Case (I)}
\label{fig: case I}
\end{figure}

Let $L$ be the line $W_1W_2$. In the figure above there are 5 configuration lines. The remaining lines $L_1, \ldots, L_6$ must form a star configuration and their intersection points $P_{ij}=L_i\cap L_j$ have to distribute in five collinear triples lying on the lines $L$, $M_1$, $M_2$, $N_1$ and $N_2$. Up to renumbering the points the collinear triples are
\begin{equation}\label{eq: cond1}
   \renewcommand{\arraystretch}{1.2}
   \begin{tabular}{|c|}
   \hline
   $P_{12}$, $P_{34}$, $P_{56}$  \\
   \hline
   $P_{13}$, $P_{25}$, $P_{46}$  \\
   \hline
   $P_{14}$, $P_{26}$, $P_{35}$  \\
   \hline
   $P_{15}$, $P_{24}$, $P_{36}$  \\
   \hline
   $P_{16}$, $P_{23}$, $P_{45}$  \\
   \hline
   \end{tabular}
   \end{equation}
   Indeed, the first column contains the points lying
   on the line $L_1$. The first row is then completed
   just by assigning numbers. The index $2$ must appear
   somewhere in the second row. Since the pairs $3,4$
   and $5,6$ cannot be distinguished in this stage (similarly
   as the particular points within these pairs), we have
   the freedom to label that point $P_{25}$. These labeling
   determines the rest of the table.

   Without loss of generality as in \eqnref{eq: abcd} and \eqnref{eq: abcd det} we may assume that
\begin{equation*}
   \begin{array}{lll}
   L_1:\; x,          & L_2:\; y,        & L_3:\; z,     \\
   L_4:\; x+y+z,      & L_5:\; ax+by+z,  & L_6:\; cx+dy+z,
   \end{array}
\end{equation*}
   with
\begin{equation*}
   a,b,c,d\in\mathbb{F}^*\;\;\mbox{ and }\; \det\left(\begin{matrix}1&1&1\\a&b&1\\c&d&1\end{matrix}\right)\neq 0.
\end{equation*}

   We can now compute coordinates of all points $P_{ij}$ for $1\leq i<j\leq 6$ and evaluate collinearity
   conditions \eqnref{eq: cond1}. This leads to the following system of equations:
\begin{equation}
\left\{\begin{array}{r} a-b-c+d=0\\-ad+a-c+d=0\\a-bc=0\\bc-d=0\\ad-a+b-d=0\end{array}\right.
\end{equation}
   This system has a solution satisfying conditions \eqnref{eq: abcd det}
   only if the ground field $\F$ has characteristic $2$. In that case we have
   $$a=d=\eps,\;\; b=c=\eps^2,$$
   with $\eps$ a solution of the equation $x^2+x+1=0$, i.e. a primitive root
   of unity of order $3$. This implies that the equations of the lines
   $L$, $N_1$, $N_2$, $M_1$, $M_2$ are (up to ordering)
   $$x+y,\;\;\; \eps x+z,\;\;\; \eps^2x+y+z,\;\;\; x+\eps^2y+z,\;\;\; \eps y+z.$$
   These lines belong all to the pencil of lines passing through the point $(1:1:\eps)$,
   a contradiction.

   Now we pass to the second case (II). In this situation we start with the following figure
\begin{figure}[H]
\centering
\begin{tikzpicture}[line cap=round,line join=round,x=1.0cm,y=1.0cm,scale=0.75]
\clip(-4.5,-0.9) rectangle (5.0,6.0);
\draw [domain=-4.5:5.0] plot(\x,{(--7.7--0.9*\x)/2.54});
\draw [domain=-4.5:5.0] plot(\x,{(--2.24-3.02*\x)/0.96});
\draw [domain=-4.5:5.0] plot(\x,{(--3.07-3.92*\x)/-1.58});
\draw [domain=-4.5:5.0] plot(\x,{(--7.02-0.56*\x)/3.85});
\draw [domain=-4.5:5.0] plot(\x,{(-2.59--3.4*\x)/0.11});
\draw [domain=-4.5:5.0] plot(\x,{(--10.04--0.36*\x)/4.22});
\begin{scriptsize}
\draw [fill=qqqqff] (-0.2,2.96) circle (2.4pt);
\draw [fill=qqqqff] (2.34,3.86) circle (2.4pt);
\draw[color=black] (-3.84,1.4) node {$N_1$};
\draw [fill=qqqqff] (0.76,-0.06) circle (2.4pt);
\draw[color=qqqqff] (1.17,0.02) node {$W_1$};
\draw[color=black] (-0.65,5.5) node {$M_1$};
\draw[color=black] (2.64,5.5) node {$M_2$};
\draw [fill=qqqqff] (-2.42,2.17) circle (2.4pt);
\draw[color=qqqqff] (-2.44,2.55) node {$W_2$};
\draw [fill=qqqqff] (1.43,1.61) circle (2.4pt);
\draw[color=black] (-3.77,2.64) node {$N_2$};
\draw [fill=qqqqff] (0.87,3.34) circle (2.4pt);
\draw[color=qqqqff] (0.62,3.62) node {$Z_3$};
\draw[color=black] (1.28,5.7) node {$M_3$};
\draw [fill=qqqqff] (1.8,2.53) circle (2.4pt);
\draw[color=qqqqff] (1.63,2.83) node {$Z_4$};
\draw[color=black] (-4.17,1.8) node {$N_3$};
\draw [fill=qqqqff] (0.82,1.7) circle (2.4pt);
\draw[color=qqqqff] (1.09,1.95) node {$Z_1$};
\draw [fill=qqqqff] (-0.01,2.38) circle (2.4pt);
\draw[color=qqqqff] (0.2,2.64) node {$Z_5$};
\draw [fill=qqqqff] (0.84,2.45) circle (2.4pt);
\draw[color=qqqqff] (1.11,2.73) node {$Z_2$};
\draw [fill=qqqqff] (0.17,1.8) circle (2.4pt);
\end{scriptsize}
\end{tikzpicture}
\caption{$ $ : Case (II)}
\label{fig: case II}
\end{figure}

   There are now $5$ remaining configuration lines, which we call as usual $L_1, \ldots, L_5$.
   They form a star configuration, i.e. there are 10 mutually distinct intersection points $P_{ij}=L_i\cap L_j$, for $1\leq i<j\leq 5$. Up to renumbering these points are distributed as follows
\begin{table}[H]
   \centering
   \renewcommand{\arraystretch}{1.2}
   \begin{tabular}{|c|c|}
   \hline
   line & points on the line \\
   \hline
   $M_1$ & $P_{12}$, $P_{34}$ \\
   \hline
   $M_2$ & $P_{15}$, $P_{23}$ \\
   \hline
   $N_3$ & $P_{13}$  \\
   \hline
   $N_1$ & $P_{14}$, $P_{25}$  \\
   \hline
   $N_2$ & $P_{24}$, $P_{35}$  \\
   \hline
   $M_3$ & $P_{45}$ \\
   \hline
   \end{tabular}
   \caption{}
   \label{tab: NiMi}
\end{table}
   The Table \ref{tab: NiMi} is filled as follows. The first two rows are filled
   just by assigning labels. They imply immediately that $Z_4\in L_4$
   and $Z_5\in L_5$. There is another triple point of the configuration
   not depicted in Figure \ref{fig: case II}. The line $L_2$ cannot
   pass through this point so that it must be $P_{13}$ (note that $L_4$ and $L_5$ intersect
   $N_3$ already in $Z_4$ and $Z_5$ respectively). The points $P_{24}$ and $P_{25}$
   must then lie one on $N_1$ and the other on $N_2$. We have selected the labeling
   in such a way, that $Z_i\in L_i$ holds for all $i=1,\ldots,5$.

   Similarly as in case (I), without loss of generality, we may assume that the equations of the lines $L_i$ are
\begin{equation*}
   \begin{array}{lll}
   L_1:\; x,          & L_2:\; y,        & L_3:\; z,     \\
   L_4:\; x+y+z,      & L_5:\; ax+by+z. &
   \end{array}
\end{equation*}
   Evaluating the conditions $Z_i\in L_i$ for $i=1,\ldots,5$
 we obtain the following system of equations
   (one condition, $Z_2\in L_2$, is satisfied automatically)
\begin{equation}
\left\{\begin{array}{r} -a^2+ab^2+ab-b^2=0\\a^2-ab^2+ab-a=0\\-ab^2+ab-a+b^2=0\\a^2-ab-a+b^2=0\end{array}\right.
\end{equation}
   This system is equivalent to the system
\begin{equation}
   \left\{\begin{array}{r} 3(a-b^2)=0\\ab-2b^2+a=0\\a^2-ab+b^2-a=0\end{array}\right.
\end{equation}
   The latter system has to be treated differently in case of characteristic of $\F$ equal either $2$ or $3$
   but all cases lead to the same solution $a=b=1$, we omit the details. This solution means $L_4=L_5$, a contradiction.

   We conclude the proof of Theorem \ref{thm: main} with an example of a configuration of
   $11$ lines with $16$ triple points. Our example is constructed over an arbitrary field $\F$
   which contains the golden section ratio. Our example is dual to the example in \cite[page 398, figure (i)]{BGS74}.

   Turning to details, let $b\in \F$ satisfy $b^2+b-1=0$ and let the configuration lines be given
   by the following equations:
\begin{equation*}
   \begin{array}{llll}
    L_1:\; x,          & L_2:\; y,        & L_3:\; z,    & L_4:\; x+y+z,\\
    L_5:\; -bx+z,        & L_6:\;bx+y+bz,    & L_7:\;y+z,   &L_8:\;b^2x+by+z,\\
    L_9:\;bx-by-z,       & L_{10}:\;-b^3x-y-bz,    & L_{11}:\;-b^2x+(1-b)y.
   \end{array}
\end{equation*}
   The configuration points are then easily computed to be:
\begin{equation*}
   \begin{array}{lll}
      P_1=(0:-1:1),              & P_2=(1:0:0), & P_3=(0:1:0),       \\
      P_4=(1:0:-1),     &P_5=(-1:b+1:-b),              & P_6=(-1:b:0),        \\
      P_7=(1:0:b),   &P_8=(0:-b:1),   &P_9=(0:-1:b), \\
      P_{10}= (1:0:-b^2), &P_{11}=(1-b:-b:b),       &P_{12}=(-1:b:-b),  \\
      P_{13}=(1:1:-1)    &P_{14}=(0:0:1),    &P_{15}=(1:1:0),\\
      P_{16}=(1:1:-2).\\
   \end{array}
\end{equation*}
   The incidence table in this case reads:
 \begin{center}
   \renewcommand{\arraystretch}{1.2}
   \begin{tabular}{|c|c|c|c|c|c|c|c|c|c|c|c|c|c|c|c|c|c|}
   \hline
   & $P_1$ & $P_2$ & $P_3$ & $P_4$ & $P_5$ & $P_6$ & $P_7$ & $P_8$ & $P_9$ & $P_{10}$ & $P_{11}$ & $P_{12}$ & $P_{13}$ & $P_{14}$ &$P_{15}$ &$P_{16}$ \\
   \hline
   $L_1$ & + & &+ & & &  &  & + &+  &  &  &  & &  +&  &  \\
   \hline
   $L_2$ & &+ &  &+  &  & &+ & &  & + &  &  & &+  &  & \\
   \hline
   $L_3$ &  & + &+ &  &  &+  &  &  & & & &  & & &+  & \\
   \hline
   $L_4$ &+  &  &  &+ &+  & &  &  & &  &  & & &  &  &+ \\
   \hline
   $L_5$ &  &  &+  &  &+ &  &+ &  &  & &  &+ &  &  &  & \\
   \hline
   $L_6$ &  &  &  &+  &  &+  &  &+ &  &  &+ & & &  &  & \\
   \hline
   $L_7$ & + &+  &  &  &  &  &  &  & &  &+  &+  &+ &  &  & \\
   \hline
   $L_8$ &  &  & &  & + & + & &  &+  &+ &  & & +&  &  & \\
   \hline
   $L_9$ & &  &  & &  &  &+  & & + & &+&  & &  & + & \\
   \hline
   $L_{10}$ & &  &  &  & & &  &+  &  &+  & & + & &  &  &+ \\
   \hline
   $L_{11}$ & & &  &  & & &  &  &  &  & &  & +&  +&  +&+ \\
   \hline
 \end{tabular}
 \end{center}
   And finally, the configuration is visualized in Figure \ref{fig: 11-16}.
   In this figure the dashed circle indicates the line at the infinity on which
   parallel lines intersect. So that
   for example parallel lines $L_9$ and $L_{11}$ intersect in point $P_{15}$.
   The line at infinity \emph{is} a configuration line.
\begin{figure}[H]
\centering
\begin{tikzpicture}[line cap=round,line join=round,x=1.0cm,y=1.0cm, scale=1.1]
\clip(-4.02,-2.6) rectangle (3.5,3.44);
\draw [domain=-4.02:9.24] plot(\x,{(--1-0*\x)/1});
\draw (1,-6.42) -- (1,3.44);
\draw (0,-6.42) -- (0,3.44);
\draw [domain=-4.02:9.24] plot(\x,{(-0-0*\x)/-1});
\draw [domain=-4.02:9.24] plot(\x,{(--0.66--1*\x)/0.65});
\draw [domain=-4.02:9.24] plot(\x,{(-0.6--1.54*\x)/1});
\draw [domain=-4.02:9.24] plot(\x,{(--0.39--0.6*\x)/-0.65});
\draw [domain=-4.02:9.24] plot(\x,{(--1-1*\x)/2.77});
\draw [domain=-4.02:9.24] plot(\x,{(-0--1.5*\x)/-1});
\draw [domain=-4.02:9.24] plot(\x,{(-0.13--0.36*\x)/-0.38});
\draw [dashed, rotate around={-160.84:(-0.03,0.27)}] (-0.03,0.27) ellipse (2.41cm and 2.29cm);
\begin{scriptsize}
\draw[color=black] (-0.22,3.25) node {$L_1$};
\draw[color=black] (-3.84,0.85) node {$L_2$};
\draw[color=black] (-0.69,2.65) node {$L_3$};
\draw[color=black] (-1.86,3.25) node {$L_4$};
\draw[color=black] (0.78,3.25) node {$L_5$};
\draw[color=black] (-2.75,3.25) node {$L_6$};
\draw[color=black] (-3.84,-0.17) node {$L_7$};
\draw[color=black] (-3.84,2.76) node {$L_8$};
\draw[color=black] (2.77,3.25) node {$L_9$};
\draw[color=black] (-3.74,1.57) node {$L_{10}$};
\draw[color=black] (1.81,3.25) node {$L_{11}$};
\fill [color=black] (0,1) circle (2.pt);
\draw[color=black] (0.39,1.2) node {$P_{14}$};
\fill [color=black] (1,0) circle (2.pt);
\draw[color=black] (1.27,0.15) node {$P_{12}$};
\fill [color=black] (1,1) circle (2.pt);
\draw[color=black] (1.375,1.2) node {$P_7$};
\fill [color=black] (0,-0.6) circle (2.pt);
\draw[color=black] (-0.26,-0.6) node {$P_9$};
\fill [color=black] (-0.66,0) circle (2.pt);
\draw[color=black] (-0.6,-0.29) node {$P_{13}$};
\fill [color=black] (0.37,0) circle (2.pt);
\draw[color=black] (0.455,-0.29) node {$P_{11}$};
\fill [color=black] (-1.77,1) circle (2.pt);
\draw[color=black] (-1.5,1.2) node {$P_{10}$};
\fill [color=black] (-0.34,0.48) circle (2.pt);
\draw[color=black] (-0.7,0.45) node {$P_{16}$};
\fill [color=black] (1,-1.5) circle (2.pt);
\draw[color=black] (1.17,-1.34) node {$P_5$};
\fill [color=black] (0,0) circle (2.pt);
\draw[color=black] (-0.17,-0.17) node {$P_1$};
\fill [color=black] (0,0.36) circle (2.pt);
\draw[color=black] (0.17,0.51) node {$P_8$};
\fill [color=black] (-0.67,1) circle (2.pt);
\draw[color=black] (-0.5,1.2) node {$P_4$};
\fill [color=black] (2.37,0.46) circle (2.pt);
\draw[color=black] (2.54,0.62) node {$P_2$};
\fill [color=black] (1.34,2.2) circle (2.pt);
\draw[color=black] (1.6,2.36) node {$P_{15}$};
\fill [color=black] (-1.88,1.66) circle (2.pt);
\draw[color=black] (-1.95,1.9) node {$P_6$};
\fill [color=black] (0.52,-1.95) circle (2.pt);
\draw[color=black] (0.59,-1.75) node {$P_3$};
\end{scriptsize}
\end{tikzpicture}
\caption{}
\label{fig: 11-16}
\end{figure}
\endproof

\section{Configurations with many triple points}\label{sec:triple configs}
   In complex algebraic geometry a point where exactly two lines meet
   is a \emph{node}. This is the simplest singularity one encounters
   and is denoted by $A_1$ in the $A$-$D$-$E$--classification of simple
   singularities of curves, see for example \cite{Arn99}. Plane curves
   (not necessarily splitting in lines) with $A_1$ singularities are
   well understood, see for example \cite{Ran86}. When exactly three lines
   meet in a point, then there is a $D_4$ singularity in that point.
   Apart from $A_1$, this is the only \emph{simple singularity} which can
   appear in an arrangement of lines. Plane curves containing $D_4$
   singularities are way less understood, see for example \cite[Section 11]{Lan03}.
   Results of this note can be considered as a step towards completing
   this picture.

   The construction we present here is directly motivated by the
   passage from the Hesse configuration to its dual.

   Let $E$ be an elliptic curve embedded as a smooth plane cubic.
   The group law on $E$ is related to the embedding by the following
   equivalent conditions
   \begin{itemize}
   \item[a)] the points $P$, $Q$ and $R$ on the curve $E$ are collinear;
   \item[b)] $P+Q+R=0$ in the group $E$.
   \end{itemize}

   Let $p$ be a prime number $\geq 3$.
   There are exactly $p^2$ mutually distinct solutions to the
   equation $pX=0$ on $E$. These solutions form a subgroup $E(p)$ of
   $p$--torsion points. Since they form a subgroup and by the above
   equivalence, any line joining two distinct points in $E(p)$
   intersects $E$ in another point which is also an element of $E(p)$.
   The tangent line to $E$ at $0$ is tangent there to order $3$, in particular
   the equation $2X=0$ on $E$
   has no non-trivial solution in $E(p)$.
   The tangent lines to $E$ at \emph{every other} point $X\in E(p)\setminus\left\{0\right\}$
   intersects $E(p)$ in some other point $Y$. This is because the equation
   $2X+Y=0$ has a unique solution in $E(p)$ for all $Y\neq 0$. In particular
   the point $X$ also lies on a line tangent to $E$ at some point $Z\in E(p)$.

   Hence, there are
   altogether $\frac{(p^2+4)(p^2-1)}{6}$ lines determined by pairs of points in
   $E(p)$. There are $\frac{p^2-1}{2}$ configuration lines
   passing through $0$ and $\frac{p^2+1}{2}$ lines passing through every other point in $E(p)$.

   Passing to the dual configuration, we obtain thus $p^2$ lines
   with $t_3(p):=\frac{(p^2-1)(p^2-2)}{6}$ triple points (and $p^2-1$ double
   points corresponding to the tangents at points $X\in E(p)\setminus\left\{0\right\}$).
   The equality in \eqnref{eq: combinatorial} guarantees that
   there are no other intersection points between the lines.
   Since $p$ is a prime, there is no rounding in \eqnref{eq: schoenheim}
   and it cannot be $p^2\equiv 5\mod (6)$,
   thus the difference $U_3(p^2)-t_3(p)=\frac{p^2-1}{3}$.

   Hence we have proved our final result.
\begin{theorem}\label{thm: main2}
   For any prime number $p\geq 3$, there exists a configuration
   of $p^2$ lines intersecting in $\frac{(p^2-2)(p^2-1)}{6}$
   triple points (and $p^2-1$ double points).
\end{theorem}
\paragraph*{\emph{Acknowledgement.}}
   This notes originated in a workshop on Arrangements of Lines held
   in Lanckorona in April 2014. We thank the Jagiellonian University in Cracow
   for financial support. We thank also Brian Harbourne and Witold Jarnicki for helpful discussions.


\bigskip \small

\bigskip
   Marcin Dumnicki, \L ucja Farnik, Halszka Tutaj-Gasi\'nska,
   Jagiellonian University, Institute of Mathematics, {\L}ojasiewicza 6, PL-30-348 Krak\'ow, Poland

\nopagebreak
   \textit{E-mail address:} \texttt{Marcin.Dumnicki@im.uj.edu.pl}

   \textit{E-mail address:} \texttt{Lucja.Farnik@im.uj.edu.pl}

   \textit{E-mail address:} \texttt{Halszka.Tutaj@im.uj.edu.pl}

\bigskip
   Agata~G\l \'owka, Magdalena~Lampa-Baczy\'nska, Grzegorz~Malara, Tomasz Szemberg, Justyna Szpond,
   Instytut Matematyki UP,
   Podchor\c a\.zych 2,
   PL-30-084 Krak\'ow, Poland

\nopagebreak


  \textit{E-mail address:} \texttt{agata.habura@wp.pl}

  \textit{E-mail address:} \texttt{lampa.baczynska@wp.pl}

   \textit{E-mail address:} \texttt{gmalara@up.krakow.pl}

   \textit{E-mail address:} \texttt{szemberg@up.krakow.pl}

   \textit{E-mail address:} \texttt{szpond@up.krakow.pl}

\bigskip

   Magdalena Lampa-Baczy\'nska, Grzegorz Malara and Tomasz Szemberg current address:
   Albert-Ludwigs-Universit\"at Freiburg,
   Mathematisches Institut, D-79104 Freiburg, Germany

\end{document}